\documentclass{article}

{}
{}
{}
\usepackage{color}
\usepackage{multirow}
\usepackage{makecell}
\usepackage{float}
\usepackage{setspace}
\usepackage{array}
\usepackage{longtable}
\usepackage{supertabular}
\usepackage{geometry}
\usepackage{fancyhdr} 
\usepackage{mathtools}
\usepackage{amssymb}
\usepackage{graphicx}
\usepackage{lineno,hyperref}

\usepackage{amsfonts}
\usepackage{comment}

\begin{document}

%\supertitle{Research Article}

\title{An Efficient MILP Formulation of Economic Dispatch with Adjustable Transformer Ratio and Phase Shifter}

\author{Bin Liu\thanks{Jibei Electric Power Dispatching and Control Center, State Grid Cooperation of China, Beijing 100054, China. Email address: eeliubin@hotmail.com.}, 
	Bingxu Zhai\thanks{Jibei Electric Power Dispatching and Control Center, State Grid Cooperation of China, Beijing 100054, China.},
 	Haibo Lan\thanks{Jibei Electric Power Dispatching and Control Center, State Grid Cooperation of China, Beijing 100054, China.}}
%\thanks{Department of Electrical Engineering, Tsinghua University, Beijing 100084, China}
%\author{\au{Bin Liu$^{1,2\corr}$}, \au{Wei Wei$^{2}$}, \au{Feng Liu$^{2}$}}
%\add{2}{Department of Electrical Engineering, Tsinghua University, Beijing 100084, China}
%\add{3}{State Grid Jibei Electric Power Dispatching and Control Center, Beijing 100053, China}
%\email{eeliubin@hotmail.com}}

\maketitle
\thispagestyle{fancy}          %更改plain状态
\fancyhead{}                      %清除以前的命令
%\lhead{This paper is a preprint of a paper submitted to The Journal of Engineering. If accepted, the
%	copy of record will be available at the IET Digital Library.}           %左上角添加
\chead{}
\rhead{}
\lfoot{}
\cfoot{}   %current page number
\rfoot{\thepage}
\renewcommand{\headrulewidth}{0pt}       %把页眉线的宽度设为零，即去掉页眉线
\renewcommand{\footrulewidth}{0pt}

\begin{abstract}
In this short paper, we study the economic dispatch with adjustable transformer ratio and phase shifter, both of which, along with the transmission line, are formulated into a generalized branch model. Resulted nonlinear parts are thereafter exactly linearized using the piecewise liner technique to make the derived ED problem computationally tractable. Numerical studies based on modified IEEE systems demonstrate the effectiveness of the proposed method to efficiency and flexibility of power system operation.  
\end{abstract}

\section{Introduction}
Economic dispatch (ED) optimizes generations of all units to achieve special purpose, e.g. economy and low carbon emission\cite{ref-101,ref-102,ref-102-1}. In traditional ED formulation based on direct current (DC) power flow, the transformer ratios (TRs) are usually approximated by 1 {\it p.u.} and the phase shifters (PSs) are not always included. Accordingly, ED is formulated as a linear programming (LP) problem, which can be efficiently solved by commercial solvers\cite{ref-103}. Although formulation of TRs and PSs are explicitly in classic optimal power flow (OPF) problem based on alternate current (AC) power flow\cite{ref-104,ref-105,ref-105-1}, the robustness of algorithms to solve AC-OPF problem is still a concern in practical power system operation. As TRs and PSs can be easily regulated on-line in practical power system operation, incorporation of such devices into traditional ED can obviously enhance the flexibility of power system, hence improving the operation efficiency. However, nonlinearity  will also be introduced into the existing linear model, making the derived ED an mixed-integer nonlinear programming (MINLP) problem, which is challenging in computation. 

In this short paper, the incorporation of TRs and PSs into traditional ED problem is studied. The introduced adjustable devices, along with transmission line, are formulated as a generalized branch model (GBM), which is then exactly linearized using piecewise linear technique (PLT). Such linearization makes the derived ED a computationally tractable mixed-integer linear programming (MILP) problem, but the much more complicated MINLP problem. Numerical studies demonstrate the effectiveness of the proposed method.

%Section II
\section{ED with Adjustable TRs and PSs}
\begin{figure}[ht]\centering\includegraphics[scale=0.8]{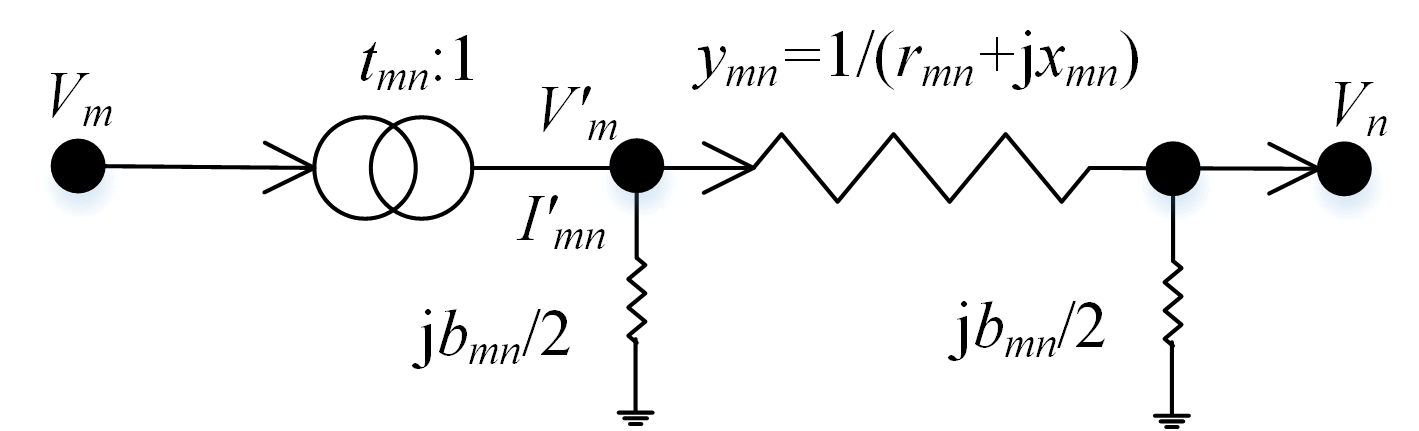}
	\caption{Generalized branch model of a single line}
	\label{fig-1}
\end{figure}

Before presenting the ED formulation with adjustable TRs and PSs, the GBM is presented based on Fig. \ref{fig-1} as follows.
%fig-1

Denoting $\mathcal{R}(\cdot)$ as the operator of returning the real part of a complex number, the active power of line $mn$ at the {\it sending} end, denoted as $P_{mn}$, can be expressed as 
\begin{equation}\label{eeq-1}\begin{split}
P_{mn}%&=\mathcal{R}\{V_m I_{mn}^*\}=\mathcal{R}\{\frac{V_m I_{mn}^{'*}}{t_{mn}} \}\\
&=\mathcal{R}\{\frac{V_m}{t_{mn}}[\frac{jV_m'b_{mn}}{2}+(\frac{V_m}{t_{mn}}-V_n) y_{mn}]^* \}
\end{split}\end{equation}
where $V$, $I$ represent the bus voltage and branch current, respectively; $y_{mn}$, $r_{mn}$, $x_{mn}$ and $b_{mn}$ are the admittance, resistance, reactance and charging capacitance of line $mn$, respectively; $t_{mn}$ is defined as $\tau_{mn} e^{j\delta_{mn}}$ with $\tau_{mn}$ and $\delta_{mn}$ representing TR and PS angle of line $mn$, respectively.

In direct-current constrained power flow (DC-PF), employing the assumptions
\begin{equation}\label{eeq-2}\begin{split}
V_m\approx e^{j\theta_m},V_n\approx e^{j\theta_n},b_{mn}=r_{mn}\approx0,\sin\theta_{mn}\approx\theta_{mn}
\end{split}\end{equation}
we have\cite{ref-103}
\begin{equation}\label{eeq-3}\begin{split}
P_{mn}&\approx\mathcal{R}\{\frac{e^{j\theta_m}}{\tau_{mn} e^{j\delta_{mn}}}[(e^{j\theta_n }-\frac{e^{j\theta_m}}{\tau_{mn} e^{j\delta_{mn}}})\frac{j}{x_{mn}} ]^* \}\\
&=\mathcal{R}\{\frac{j}{\tau_{mn} x_{mn}}[\frac{1}{\tau_{mn}} -e^{j(\theta_m-\theta_n-\delta_{mn} )}]\}\\
&=\mathcal{R}\{\frac{j}{\tau_{mn} x_{mn}} [\frac{1}{\tau_{mn}}-\cos⁡{(\theta_m-\theta_n-\delta_{mn})}\\
&-j\sin{⁡(\theta_m-\theta_n-\delta_{mn} )} ]\}\\
&\approx\frac{\theta_m-\theta_n-\delta_{mn}}{\tau_{mn} x_{mn} }
\end{split}\end{equation}

Based on \eqref{eeq-3}, the ED problem with multiple periods can be formulated as \eqref{eeq-4}-\eqref{eeq-7_4}, where the objective, which is to minimize the overall operation cost is presented as \eqref{eeq-4} in compact form and the constraints of generation capacities, ramping rates and reserve requirements are presented in compact form as \eqref{eeq-5}. The details of \eqref{eeq-4} and \eqref{eeq-5} in the ED model can be found in \cite{ref-103} and is omitted here for simplicity. Power balance and transmission line capacity constraints in the formulation are presented in as \eqref{eeq-6} and \eqref{eeq-7}. The total allowable times of adjusting TR and PS are limited by \eqref{eeq-7_1}-\eqref{eeq-7_3} and the adjusting steps of TR and PS are constrained by \eqref{eeq-7_4}-\eqref{eeq-7_5}.
\begin{eqnarray}
\label{eeq-4}&\min_{p\in\mathcal{S}}{f(p)}\\
\label{eeq-5}&s.t.~\mathcal{S}=\{p|\mathbf{A}p+\mathbf{B}d\le a\}\\
\label{eeq-6}&p_{m,h}-d_{m,h}=\sum\nolimits_{n\in\mathcal{N}_m}{P_{mn,h}}~~\forall m,\forall h\\
\label{eeq-7}&-\overline{L}_{mn}\le P_{mn}=\frac{\theta_{m,h}-\theta_{n,h}-\delta_{mn,h}}{x_{mn}\tau_{mn,h}}\le\overline{L}_{mn}~\forall mn,\forall h\\
\label{eeq-7_1}&|\tau_{mn,h}-\tau_{mn,h-1}|\le I^{\tau}_{mn,h}(\overline\tau_{mn}-\underline\tau_{mn})~\forall mn,\forall h\\
\label{eeq-7_2}&|\delta_{mn,h}-\delta_{mn,h-1}|\le I^{\delta}_{mn,h}(\overline\delta_{mn}-\underline\delta_{mn})~\forall mn,\forall h\\
\label{eeq-7_3}&\sum\nolimits_{h}{I^{\tau}_{mn,h}}\le n_{mn}^{\tau},~\sum\nolimits_{h}{I^{\delta}_{mn,h}}\le n_{mn}^{\delta}~\forall mn\\
\label{eeq-7_4}&|\tau_{mn,h}-\tau_{mn,h-1}|\le \Delta\tau_{mn}~\forall mn,\forall h\\
\label{eeq-7_5}&|\delta_{mn,h}-\delta_{mn,h-1}|\le \Delta\delta_{mn}~\forall mn,\forall h
\end{eqnarray}
where $m, h$ and $mn$ are indices of unit, dispatch interval and branch; $f(p)$ represents the fuel cost of thermal units and can be reformulated as linear functions\cite{ref-106}; $p$ and $d$ represent power generation and nodal power load; $\mathcal{N}_m$ represents the neighboring bus set of $m$; $\overline{L}_{mn}$ is the transmission capacity of $mn$; $\underline{\cdot}$ and $\overline{\cdot}$ represent lower and upper bounds of corresponding variable or parameter; $I^{\tau}$ and $I^{\delta}$ are intermediate binary variables; $n^{\tau}_{mn}$ and $n^{\delta}_{mn}$ are total allowable times of adjusting TR and PS during the whole dispatching period for line $mn$; $\Delta{\tau}_{mn}$ and $\Delta{\delta}_{mn}$ are maximum step size of adjusting TR and PS for line $mn$, respectively; $\mathbf A, \mathbf B,a$ are constant parameters in appropriate dimensions.% and $P_{mn}$ is expressed as \eqref{eeq-3}.% $n_{mn}^{\tau},n_{mn}^{\delta}$ are total allowable adjusting times of TR and PS for $mn$, respectively; $\Delta\tau_{mn},\Delta\tau_{mn}$ are change rates of TR and PS for $mn$, respectively;

Noting that the nonlinear constraints \eqref{eeq-7_1}, \eqref{eeq-7_2}, \eqref{eeq-7_4} and \eqref{eeq-7_5} can be easily reformulated to linear forms by eliminating the absolute operator, the remaining nonlinearity is due to the formulation of $P_{mn}$ in \eqref{eeq-7}, which makes the problem computationally intractable.

%Section III
\section{Exactly Linearizing GBM in ED}
The non-convex and nonlinear part, i.e. $P_{mn}$, should be reformulated to make the problem computationally tractable, which is illustrated as follows taking branch $mn$ as an example. For branch $mn$, let 
\begin{equation}\label{eeq-8}\begin{split}
\theta_m\in[\underline\theta_m,\overline\theta_m];\theta_n\in[\underline\theta_n,\overline\theta_n];\delta_{mn}\in[\underline\delta_{mn},\overline\delta_{mn}]\\
\tau_{mn}\in\mathcal{W}_{mn}=\{\omega_{mn,1},\cdots,\omega_{mn,K_{mn}}\}
\end{split}\end{equation}
where $\mathcal{W}_{mn}$ is the discrete feasible region of $\tau_{mn}$ with $K_{mn}$ representing the total number of strategies in it. 

Introducing $6K_{mn}$ continuous variables $z^{\theta_m}_{mn,i,j}, z^{\theta_n}_{mn,i,j}$ and $z^{\delta_{mn}}_{mn,i,j}(\forall i\in\{1,\cdots,K_{mn}\},\forall j\in\{1,2\})$ and $K_{mn}-1$ binary variables $y_{mn,k}(\forall k\in\{1,\cdots,K_{mn}-1\})$, $P_{mn}$ according to PLT \cite{ref-107,ref-108} can be linearized as 
\begin{equation}\label{eeq-9}\begin{split}
P_{mn}=Y_{mn}(\theta_m)-Y_{mn}(\theta_n)-Y_{mn}(\delta_{mn})\\
\sum\nolimits_k{y_{mn,k}}=1
\end{split}\end{equation}
where $Y_{mn}(\alpha)$ with $\alpha\in[\underline\alpha,\overline\alpha]$ is defined as
\begin{equation}\label{eeq-10}\begin{split}
Y_{mn}(\alpha)=\sum\nolimits_i{(z^{\alpha}_{mn,i,2}\overline\alpha+z^{\alpha}_{mn,i,1}\underline\alpha)/(\omega_{mn,i}}x_{mn})\\
\tiny{s.t.~~~~}\alpha=\sum\nolimits_i{(z^\alpha_{mn,i,2}\overline\alpha+z^\alpha_{mn,i,1}\underline\alpha)}~~\forall j\\
z^\alpha_{mn,i,j}\ge0~~\forall i,\forall j\\
\tau_{mn}=\sum\nolimits_i^{}{(z^{\alpha}_{mn,i,2}+z^{\alpha}_{mn,i,1})\omega_{mn,i}}\\
\sum\nolimits_i^{}{\sum\nolimits_j{z^\alpha_{mn,i,j}}}=1\\
z^\alpha_{mn,1,j}\le y_{mn,1},~z^\alpha_{mn,K_{mn},j}\le y_{mn,K_{mn}-1}~\forall j\\
z^\alpha_{mn,l,j}\le y_{mn,l-1}+y_{mn,l}~\forall l\in\{2,\cdots,K_{mn}-1\},\forall j
\end{split}\end{equation}

It is noteworthy that the linearization is \textbf{exact} due to the fact that $\tau_{mn}$ is a discrete variable belonging to a finite set. Actually, for any variable $z$ that can be expressed as  $z=x^\lambda y$, where $\lambda$ is a known constant, $x$ is a discrete variable belonging to an finite set and $y$ is a continuous variable belonging to a box set, $z$ can be \textbf{exactly} linearized by PLT. The proof can be found in \cite{ref-108} and is omitted here for simplicity.

Correspondingly, ED with GBM is reformulated as a MILP problem, i.e. \eqref{eeq-4}-\eqref{eeq-9} and mixed-integer linear constraints affiliated to $Y_{mn}(\theta_m)$, $Y_{mn}(\theta_n)$ and $Y_{mn}(\delta_{mn})$, which can be solved by commercial solvers efficiently. It is also noteworthy that the employment of PLT is necessary when a branch contains an adjustable TR. However, for branches with fixed TR, either with or without any PS, the employment of PLT is unnecessary as $P_{mn}$ is itself linearly formulated.  

Besides, considering that TRs may be employed to regulate bus voltages in real-time operation, $\mathcal{W}_{mn}(\forall mn)$ can be compressed to assure that bus voltage requirements can be satisfied under any realization of $\tau_{mn}\in\mathcal{W}_{mn}(\forall mn)$, either by the dispatcher's experiences or post-simulation of the obtained operation strategy. However, this work is beyond the scope of this paper and will be studied in the future. Besides, the problem that the proposed model may give a physically infeasible solution will not be discussed in this paper as it is a concern in all DC-PF based operation problems. %Such a consideration can be realized by checking the feasibility of the AC power flow with derived optimal transformer ratios, which is beyond the scope of this letter, thus will not be discussed.  

%Scetion IV
\section{Numerical Studies}
The proposed method are studied in this section based on modified IEEE 6-bus, 39-bus and 118-bus systems, the data of which can be found in \cite{ref-103}. The locations of TRs and PSs are briefly presented in Table \ref{tab-0} with $\Delta\tau_{mn}=0.01$, $\Delta\delta_{mn}=3^0$, $\tau_{mn}\in\{0.98,0.99,\cdots,1.02\}$, $\delta_{mn}\in[-15^0,15^0]$ and $n^\tau_{mn}=n^\delta_{mn}=8$ for all $mn$ in all studied systems.

Besides, we set $\overline L_{23}=\overline L^0_{23}\times20\%$ for 39-bus system and $\overline L_{35}=\overline L^0_{35}\times50\%$ for 118-bus system, where $\overline L^0_{23}$ and $\overline L^0_{35}$ are original transmission capacities presented in \cite{ref-103}. Practically, such modification will be required by the Electric Power Dispatching and Control Center in all power companies in China under emergent situations such as overheating problems of transmission lines or switches, stability considerations etc, which need reduce the currents running through these devices to lower the risk of disconnecting them from the power system. All models are solved by CPLEX 12.6 with termination gap $0.01\%$ on a Thinkpad equipped with Intel(R) i7-3520M-2.9GHz and 8GB RAM. For comparison, ED without and with adjustable TRs and PSs are denoted by $ED_0$ and $ED_1$, respectively. All computation results are showed in Table \ref{tab-1}.% and Table \ref{tab-2}. 
%tab-0
\begin{table}[ht]\renewcommand\arraystretch{1.2}
	\caption{Locations of Adjustable TRs and PSs}
	\begin{tabular}{>{\scriptsize}c|>{\scriptsize}c|>{\scriptsize}c}
		\hline\hline
		Systems&TR Locations& PS Locations\\
		\hline
		6-bus&2, \textbf5&\textbf5, 7\\
		\hline
		%39-bus&\makecell{\textbf5, 14, \textbf{20}, 21, 22,\\32, \textbf{33}, 34, 37, 39, 41, 46}&\makecell{\textbf5, 11, \textbf{20}, 27, \textbf{33}, 43}\\
		%\hline
		39-bus&\makecell{\textbf{21}, \textbf{22}, \textbf{32}}&\makecell{5, 11, \textbf{21}, \textbf{22}, 27, \textbf{32}, 37, 44}\\
		\hline
		118-bus&\makecell{8, \textbf{32}, 36, \textbf{51},\\ \textbf{93}, 95, \textbf{102}, 107, \textbf{127}}&\makecell{24, 29, \textbf{32}, 38, \textbf{51},\\90, \textbf{93}, \textbf{102}, 105, 125, \textbf{127}}\\
		\hline\hline
	\end{tabular}
	\centering
	\label{tab-0}
\end{table}
\begin{table}[ht]\renewcommand\arraystretch{1.2}
	\caption{Computation Results ($ED_0$//$ED_1$)}
	\begin{tabular}{>{\scriptsize}c>{\scriptsize}c>{\scriptsize}c>{\scriptsize}c}
		\hline\hline
		Systems& Cost (\$)&\makecell{Time (s)}&\makecell{Cost Reduction (\%)}\\
		\hline
		{6-bus}&76687.4//72064.5&0.31//1.56&6.03\\
		\hline
		{39-bus}&479065.6//470492.1&	0.25//14.67&1.79\\
		\hline
		{118-bus}&infeasible//1838310.2&---//1255.49&---\\
		\hline\hline
	\end{tabular}
	\centering
	\label{tab-1}
\end{table}
%tab-1end

For 6-bus system, the operation cost is reduced by 6.03\%, demonstrating that incorporating adjustable TRs and PSs can bring effective economic benefit to power system operation. 
%fig-2
\begin{figure}[ht]\centering\includegraphics[scale=0.24]{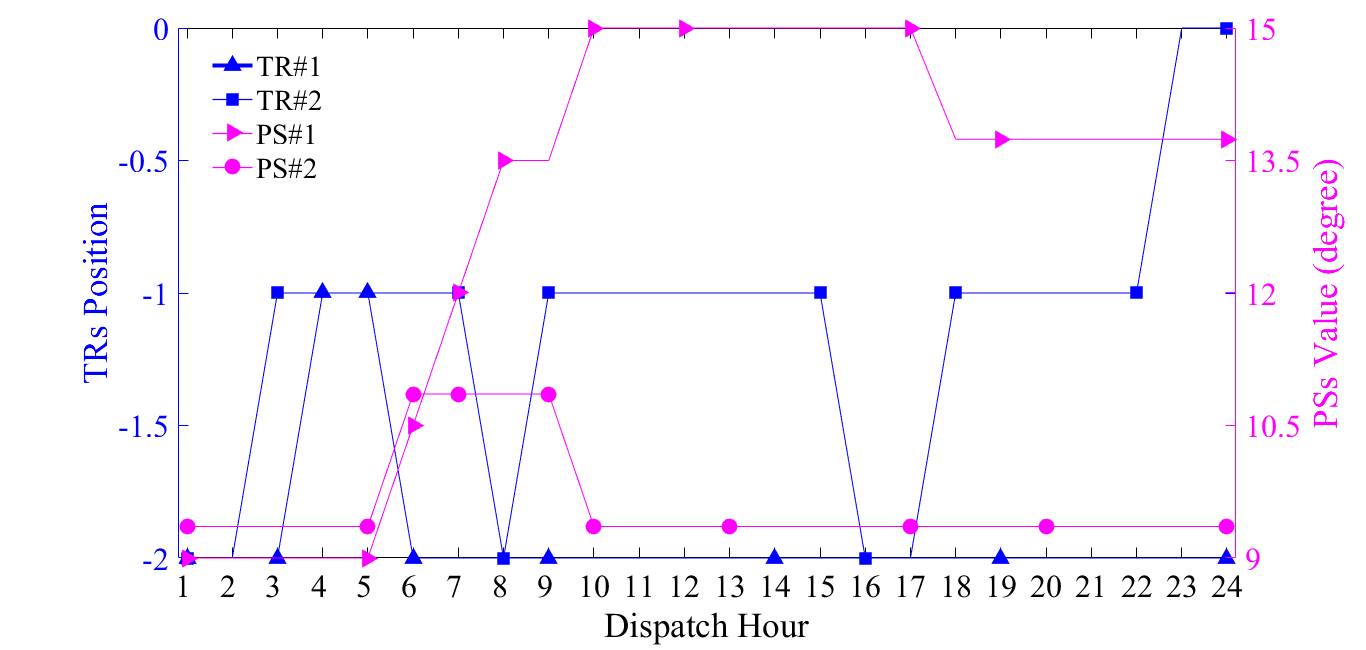}
	\caption{Optimal strategies of TRs and PSs for 6-bus system}
	\label{fig-2}
\end{figure}
%fig-3
\begin{figure}[ht]\centering\includegraphics[scale=0.24]{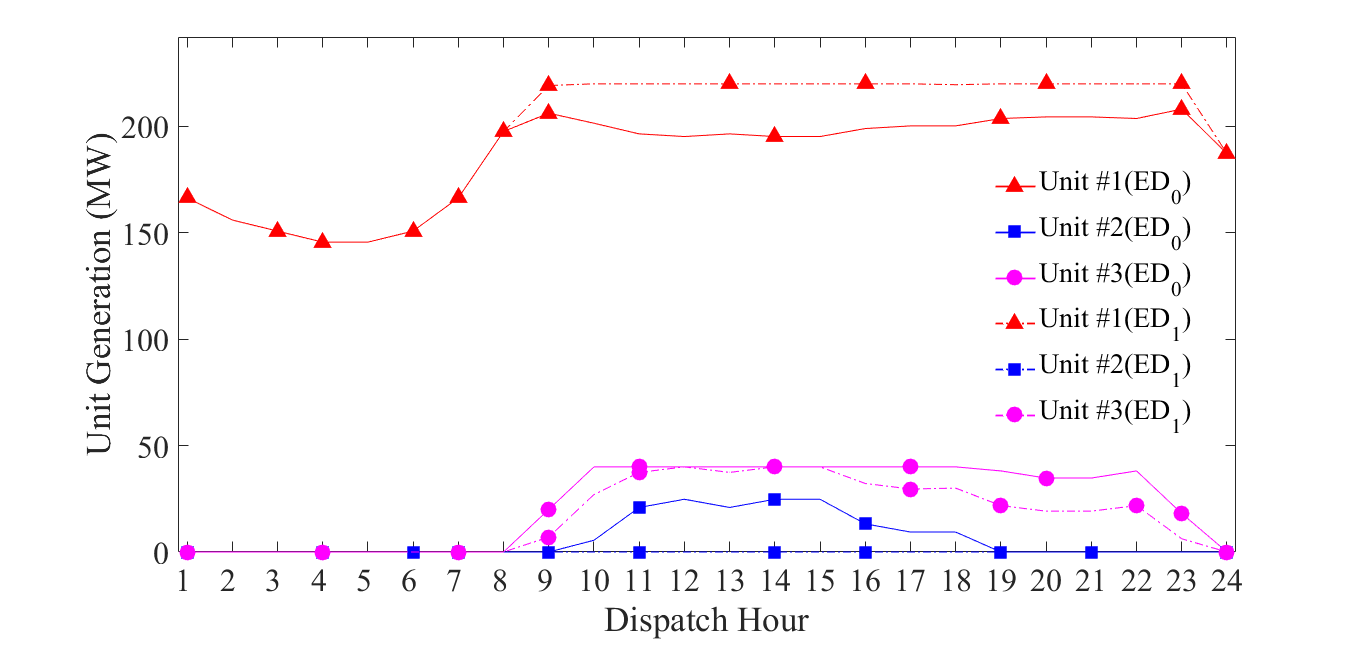}
	\caption{Optimal unit generations of $ED_0$ and $ED_1$ for 6-bus system}
	\label{fig-3}
\end{figure}

The optimal operation strategies of TRs and PSs and optimal unit generations in all dispatch hours with/without adjustable TRs and PSs are presented in Fig. \ref{fig-2} and Fig.\ref{fig-3}. Regarding the computation results presented in these two figures, more electric power will be provided by cheaper units, i.e. unit 1 and unit 3 when adjustable TRs and PSs are operated at optimal strategies. This explains how the economic benefit is brought by the introduced adjustable devices.

For 39-bus system, the operation cost is reduced by 1.79\%. However, for 118-bus system, a feasible solution can be found only if adjustable TRs and PSs are considered due to limitations of transmission capabilities of line 23 and 35, which shows that the flexibility of power system is enhanced. It is also interesting to note that before modifying the transmission line capacities, i.e. $\overline L_{23}=\overline L^0_{23}$ for 39-bus system and $\overline L_{35}=\overline L^0_{35}$ for 118-bus system, optimal unit generations and operation costs with adjustable devices are the same as that without such devices due to the fact that no transmission lines are congested. This means the proposed method may play an more important role for power systems with heavily loaded transmission lines.

%Section V
\section{Conclusions}
In this short paper, an efficiently MILP formulation of ED with adjustable TRs and PSs are presented and studied. Numerical studies based on modified IEEE systems demonstrate the effectiveness of enhancing power system operation efficiency or flexibility brought by adjustable TRs and PSs. 

% Non-BibTeX users please use

%\appendix
\end{document}